\documentclass{amsart}
\usepackage[dvips]{color}
\usepackage{graphicx}
\usepackage{epsfig}
\usepackage{tikz}
\usepackage{enumerate}
\usepackage{enumitem}
\usepackage{color}

\usepackage{latexsym}
\usepackage{amssymb}
\usepackage{amsmath}
\usepackage{amsthm}
\usepackage{amscd}
\usepackage{latexsym}
\usepackage{amssymb}
\usepackage{amsmath}
\usepackage{amsthm}
\usepackage{amscd}
\theoremstyle{plain}
\newtheorem{thm}{Theorem}[section]
\newtheorem{lem}[thm]{Lemma}

\theoremstyle{definition}
\newtheorem*{definition*}{Definition}

\newtheorem*{notation*}{Notation}

\newtheorem{cor}[thm]{Corollary}

\newcommand{\future}[1]{{}}

\newcommand{\diam}{\operatorname{diam}}%Need to add to main file
\newcommand{\cR}{{\mathcal R}}%Need to add to main file

\title{Self joinings of rigid rank one transformations arise as strong operator topology limits of convex combinations of powers} 

\author{Jon Chaika}
\thanks{The research of J. Chaika was supported in part by NSF grants DMS-135500 and DMS- 1452762, the Sloan foundation, Poincar\'e chair, Warnock chair}
\address{Department of Math
155 South 1400 East, JWB 233
Salt Lake City, UT 84112}
\email{chaika@math.utah.edu}

\begin{document}

\maketitle

The following result is a straightforward modification of \cite[Section 2]{CE} of the first named author and A. Eskin, which is included for ease of future reference. What is below is an edited version of that section due to some added (but straightforward to address) difficulties to make the assumptions be conjugacy invariant and so hold for a residual set of measure preserving transformations. For connections to the work of others see that paper. 

Let $([0,1],\mathcal{M},\lambda, T)$, be an ergodic invertible transformation. We say it is
\emph{rigid rank 1}  if there exist numbers $n_j$ and measurable sets $A_j$ such that 
%\marginpar{\bkc{in the main file, we defined rigid rank 1 by cylinders -- we could rephrase this in those terms?}}
\begin{enumerate}
\item\label{cond:big tower} $\underset{j \to \infty}{\lim}\, \lambda(\bigcup_{i=0}^{n_j-1}T^iA_j)=1$,
\item The sets $A_j,...,T^{n_j-1}A_j$ are pairwise disjoint.
\item $\underset{j \to \infty}{\lim}\, \frac{\lambda(T^{n_j}A_j \cap A_j)}{\lambda(A_j)}= 1$

\item\label{cond:bunched} For all $\varepsilon>0$ there exist metric balls of diameter at most $\varepsilon$, $B^{(j)}_0,...B^{(j)}_{n_j-1}\subset [0,1]$ such that 
$$\underset{j \to \infty}{\lim}\, \sum_{i=0}^{n_j-1}\lambda(T^iA_j\setminus B^{(j)}_i)=0.$$
\end{enumerate}
%This is a condition saying that our transformation is well
%approximated by periodic transformations. A similar condition,
%admiting \emph{cyclic approximation by periodic transformations} was
%considered in \cite{KS}.

Let
\begin{equation}\label{eq:def:Rk}
\cR_k=\bigcup_{i=0}^{n_k-1}T^iA_k,
\end{equation}
\begin{equation}
\label{eq:def:hatRk}
\hat{\cR}_k = \bigcup_{i=0}^{n_k-1} T^i( I_k \cap T^{-n_k} A_k \cap
T^{n_k} A_k),
\end{equation}
\begin{displaymath}
\tilde{\cR}_k = \bigcup_{i=0}^{n_k-1} T^i( A_k \cap T^{-n_k} A_k \cap
T^{-2n_k} A_k \cap T^{n_k} A_k \cap T^{2 n_k} A_k),
\end{displaymath}
Then, 
$\cR_k$ is the Rokhlin tower over $A_k$,  $\hat{\cR}_k$ is the Rokhlin
tower over $A_k \cap T^{-n_k}A_k\cap T^{n_k}A_k$, and $\tilde{\cR}_k$ is the Rokhlin tower over 
$\bigcap_{i=-2}^2T^{in_k}A_k$. 
We have
\begin{equation}
\label{eq:property:hatRk}
\hat{\cR}_k\supset \{x:T^ix \in \cR_k \text{ for all }-n_k< i<n_k\},
\end{equation}
and
\begin{equation}
\label{eq:property:tildeRk}
\tilde{\cR}_k\supset \{x:T^ix \in \cR_k \text{ for all }-2n_k< i< 2n_k\}. 
\end{equation}
Heuristically one can think of $\cR_k$ as the set of points we can control. $\hat{\cR}_k$ and $\tilde{\cR}_k$ let us control the points for long orbit segments, which is necessary for some of our arguments. 
\begin{lem}\label{lemma:srank est} $\underset{k \to \infty}{\lim}\lambda(\tilde{\cR}_k)=1=\underset{k \to \infty}{\lim}\lambda(\cR_k)=\underset{k \to \infty}{\lim}\lambda(\hat{\cR}_k)$. 
\end{lem}
\begin{proof}By the first condition in the definition of rigid rank 1
 we have $\underset{k \to \infty}{\lim}
\lambda(\cR_k)=1$. By (\ref{eq:def:hatRk}),
\begin{displaymath}
\lambda(\hat{\cR}_k)\geq\lambda(\cR_k)-n_k\lambda(A_k\setminus
(T^{n_k}A_k \cup T^{-n_k}A_k))
\geq\lambda(\cR_k)-2n_k\lambda(A_k\setminus T^{n_k}A_k),
\end{displaymath}
and thus by the fourth condition of the definition of rigid rank 1 by
intervals, $\underset{k \to \infty}{\lim}\lambda(\hat{\cR}_k) \to 1$. 
% By the fourth condition we have that 
% $\underset{k \to \infty}{\lim}\, \lambda(R_k \Delta \cup_{i=0}^{n_k-1}T^{i }(I_k \cap T^{n_k}I_k \cap T^{-n_k}I_k))=0 $ and so $\underset{k \to \infty}{\lim} \lambda(\hat{R}_k)=1$ because $\hat{R}_k= \cup_{i=0}^{n_k-1}T^{i }(I_k \cap T^{n_k}I_k \cap T^{-n_k}I_k)$.  
Similarly, $\underset{k \to \infty}{\lim}
\lambda(\tilde{R}_k)=1$. 
\end{proof}

%\bold{Shifted power joinings.}
%Let $(X,T,\mu)$ and $(Y,S,\nu)$ be measure preserving dynamical
%systems. Recall that a joining of $(X,T,\mu)$ and $(Y,T,\nu)$ is a
%$T \cross S$-invariant measure $\sigma$ on $X \cross Y$ which projects
%to $\mu$ under the map $X \cross Y \to X$ and to $\nu$ under the map
%$X \cross Y \to Y$. 

%\begin{defin}[Shifted Power Joining]
%Let $(X,T,\mu)$ be a measure preserving dynamical
%system. A self-joining of $(X,T,\mu)$ %with $(X,T,\mu)$
%that gives full measure to $\{(x,T^a x)\}$ for some $a\in \mathbb{Z}$ with
%$a \ne 0$ is called a \emph{shifted power joining}. 
%\end{defin}
%These have also been called \emph{off diagonal joinings}.

%Let $\iota:[0,1] \to [0,1]$ by $x \to (x,x)$. Let
%$\mu=\iota_*\lambda$. Shifted power joinings have the form $(id \times
%T^a)_* \mu$ for some $a \in \mathbb{Z} \setminus \{0\}$.   

\textbf{The operator }$A_\sigma$ \textbf{ and convergence in the strong operator topology.}
Let $\sigma$ be a self-joining of $(T,\lambda)$.  %\marginpar{\jcc{Should we change notation for the $L^2$ function $A_\sigma f$ evaluated at $x$?}}
 Let $\sigma_x$ be the corresponding measure on $[0,1]$ coming from disintegrating $\sigma$ along projection onto the first coordinate. Note this is a slight abuse, as we are identifying the measures on the fibers $\{x\}\times [0,1]$ with measure on $[0,1]$. 
Define $A_\sigma: L^2(\lambda) \to L^2(\lambda)$ by $A_\sigma(f)(x)=\int f
d\sigma_x$. %\marginpar{\jcc{Delete definition of strong operator topology?}}

Recall that one calls the \emph{strong operator topology}
the topology of pointwise convergence on $L^2(\lambda)$. That is
$A_1,...$ converges to $A_\infty$ in the strong operator topology if
and only if $\underset{ i \to \infty}{\lim}\|A_if-A_{\infty}f\|_2=0$ for all $f \in L^2(\lambda)$. 
\begin{thm} 
\label{theorem:SOT close}
Assume $([0,1],T,\lambda)$ is rigid
rank 1  and $\sigma$ is a self-joining of
$([0,1],T,\lambda)$. Then $A_\sigma$ is the strong operator topology
(SOT) limit of linear combinations, with non-negative coefficients, of powers of $U_T$, where $U_T:
L^2([0,1],\lambda) \to L^2([0,1],\lambda)$ denotes the Koopman
operator $U_T(f) = f \circ T$. 
\end{thm}
Given $n \in \mathbb{Z}$, we obtain a self-joining of $([0,1],T,\lambda)$ carried on $\{(x,T^nx)\}$, $J(n)$ defined by $\int_{X\times X} fdJ(n)=\int_X f((x,T^nx))d\mu$. We call this an \emph{off-diagonal joining}.
\begin{cor}\label{cor:WOT close} (J. King \cite{flat stacks}) Any self-joining of a rigid rank 1  transformation is a weak-* limit of linear combinations, with non-negative coefficients, of off diagonal joinings.
\end{cor}
%These results (or very closely related results) were established earlier by J. King \cite{king} using a different proof. In fact he shows that if the joining in Corollary \ref{cor:WOT close} is ergodic then there is no need to take a linear combination. See also \cite[Theorem 7.1]{JRR}. There is an open question of whether this result is true for general rank 1 systems \cite[Page 382]{King r1}. Ryzhikov has a series of results in this direction, see for example \cite{ryzh1} and \cite{ryzh2}. 
%Note that these are a partial answer to a much harder question of King 
%who asked if for rank 1 transformations shifted power joinings were dense in ergodic joinings \cite[Page 382]{King r1}.
\subsection{Proof of Theorem \ref{theorem:SOT close}}
\begin{lem}
For each $0\leq j<n_k$ we have 
\begin{equation}
\label{eq:lemma1:4}
n_k\int_{T^jA_k}\sigma_x (\cR_k^c)d \lambda(x)\leq \lambda (\hat{\cR}_k^c).
\end{equation}
\end{lem}
\textbf{Remark.} Note that $n_k$ is roughly $\lambda(T^jA_k)^{-1}$. 
\begin{proof}
Suppose $0 \leq j < n_k$, and suppose $x \in T^j A_k$. 
From (\ref{eq:property:tildeRk}) we have
$T^{i}\cR_k^c \subset \hat{\cR}_k^c$ for all $-n_k<i<n_k$. 
We claim that
\begin{equation}
\label{eq:sigmax:Rkc}
\sigma_x (\cR_k^c)\leq \sigma_{T^{\ell} x}(\hat{\cR}_k^c) \qquad \text{for all $-n_k < \ell<n_k$.} 
\end{equation}
Indeed,
$\sigma_{x}(\cR^c_k)=\sigma_{T^\ell x}(T^{\ell}\cR^c_k)\leq
\sigma_{T^\ell x}(\hat{\cR}_k^c)$, proving (\ref{eq:sigmax:Rkc}). 
Integrating (\ref{eq:sigmax:Rkc}) we get
\begin{equation}
\label{eq:int:TjIk:sigma:y}
\int_{T^jA_k}\sigma_y(\cR_k^c)d\lambda(y)\leq
\int_{T^{j+\ell}A_k}\sigma_z(\hat{\cR}_k^c) d\lambda(z) \quad\text{for all
$-n_k <\ell<n_k$.}
\end{equation}
Since we can choose $\ell$ in (\ref{eq:int:TjIk:sigma:y}) so that
$j+\ell$ takes any value in $[0,n_k-1]\cap \mathbb{Z}$, we get
\begin{equation}
\label{eq:int:Tj:Ik:min}
\int_{T^jA_k}\sigma_y(\cR_k^c)d\lambda(y)\leq \min_{0 \le i < n_k}
\int_{T^i A_k}\sigma_z(\hat{\cR}_k^c) d\lambda(z).
\end{equation}
Now
$$\sum_{i=0}^{n_k-1}\int_{T^iA_k}\sigma_y(\hat{\cR}_k^c)d\lambda(y)\leq
\int_{[0,1]}\sigma_y(\hat{\cR}_k^c)d\lambda(y)=
\lambda(\hat{\cR}_k^c),$$ 
where the last estimate
uses that $\sigma$ has projections $\lambda$. So we obtain 
\begin{equation}
\label{eq:min:small}
\min\limits_{0\leq i<n_k} \int_{T^iA_k}\sigma_x
(\hat{\cR}_k^c)d\lambda(x)\leq  \frac 1 {n_k} \lambda(\hat{\cR}_k^c).
\end{equation}
Now the estimate (\ref{eq:lemma1:4}) follows from
(\ref{eq:int:Tj:Ik:min}) and (\ref{eq:min:small}). 
\end{proof} 

We want to guess coefficients $c_j$ so that $A_\sigma$ is close to $\sum_{j=0}^{n_k-1}c_jU_T^j$. The next lemma comes up with a candidate pointwise version. Theorem \ref{theorem:SOT close} and Corollary \ref{cor:WOT close} will follow because by Egoroff's theorem this choice is almost constant on most of the $T^\ell I_k$ 
and the lemma after this (Lemma~\ref{lemma:other indices}), which shows that they are almost $T$ invariant. 
\begin{lem}\label{lemma:A close} Let $x \in \hat{\cR}_k \cap T^jA_k $ where $0\leq j<n_k$.  Define $c_i(x)=\sigma_x(T^{a_i}A_k\cap {\cR}_k)$ where $0\leq a_i<n_k$ and $i+j \equiv a_i \,( \text{mod }n_k)$.
 For all 1-Lipschitz $f$ we have 
 \begin{multline*}\left|A_\sigma f(x)-\sum_{i=0}^{n_k-1}c_i(x)f(T^ix)\right|\leq \varepsilon+\|f\|_{\sup}\sigma_x({\cR}_k^c) +\\
 \|f\|_{\sup}\sigma_x\big(\cup_{i=0}^{n_k-1}(T^iA_k \setminus B_i^{(k)})\big)+\|f\|_{\sup}\sum_{i:T^ix\notin B_{a_i}^{(k)}}\sigma_x(T^iA_k).
 \end{multline*}
\end{lem}
% Morally $c_j(x)$ is the $\sigma_x$ measure of the level in $R_k$ that is $j$ levels above the level $x$ is on. Because $j+\ell$ can be bigger than $n_k$ the definition is slightly more complicated.  Note that the $c_j(x)$ are non-negative. 
 \begin{proof}First observe that 
 \begin{equation}\label{eq:tower approx}|A_\sigma f(x)-\sum_{i=0}^{n_k-1}\int _{T^iA_k}f d\sigma_x|\leq \|f\|_{\sup}\sigma_x(\cR_k^c)
 \end{equation}
 Now if $f(T^\ell x) \in B_i^{(k)}$ we have 
 \begin{equation*}|\int_{T^iA_k}fd\sigma_x-f(T^\ell x)\sigma_x(T^iA_k)|\leq \epsilon\|f\|_{Lip}\sigma_x(B_i^{(k)})+\|f\|_{\sup}\sigma_x(T^iA_k\setminus B_i^{(k)}).
 \end{equation*}
By applying the above estimate if $T^ix\in B_{a_i}^{(k)}$ and estimating trivially if it isn't we obtain 
\begin{multline*}
|\sum_{i=0}^{n_k-1}\int_{T^iA_k}f d\sigma_x-\sum_{i=0}^{n_k-1}f(T^ix)\sigma_x(T^{a_i}A_k)|= |\sum_{i=0}^{n_k-1}\int_{T^iA_k}f d\sigma_x-\sum_{i=0}^{n_k-1}c_{a_i}f(T^ix)|\\
\leq \sum_{i=0}^{n_k-1}\|f\|_{\sup}\sigma_x(T^iA_k\setminus B^{(k)}_i)+\|f\|_{\sup}\sum_{i:T^ix\notin B_{a_i}^{(k)}}\sigma_x(T^iA_k)+\epsilon\|f\|_{Lip}
\end{multline*}
Combining this with \eqref{eq:tower approx} gives the lemma. 
 \end{proof}

\begin{lem}
\label{lemma:other indices}
Suppose $0 \le \ell < n_k$. 
If $x\in T^\ell A_k$ and $-\ell\leq i<n_k-\ell$ then 
$$\sum_{j=0}^{n_k-1}|c_j(x)-c_j(T^ix)|\leq 2\sigma_x(\tilde{R}_k^c).$$
\end{lem}
\begin{proof}
Suppose $0 \le \ell < n_k$, $0 \leq j < n_k$, and $-\ell \le i <
n_k-\ell$.  
First note that  if $0 \le m < n_k$ and $z \in T^m A_k =T^mA_k\cap {R}_k$ then by
(\ref{eq:def:Rk}), we have $T^s z \in T^{m+s}A_k \cap
{R}_k$ for all $-m\leq s<n_k-m$. Thus, if $j + \ell <
n_k$ and $i+j+\ell < n_k$, we have
$$\sigma_{T^ix}(T^{i+j+\ell}A_k \cap {R}_k)=\sigma_{x}(T^{j+\ell}A_k \cap T^{-i}{R}_k)=\sigma_x(T^{j+\ell}A_k \cap {R}_k).$$
This gives $c_j(x)=c_j(T^ix)$ if $j+\ell <n_k$ and $i+j+\ell<n_k$. 
By similar reasoning we have that $c_j(x)=c_j(T^ix)$ if $j+\ell\ge n_k$ and $i+j+\ell\geq n_k$. 

Now lets assume that $j+\ell<n_k$ and $i+j+\ell\geq n_k$. Then,
\begin{equation}
\label{eq:cjTix}
c_{j}(T^ix)=\sigma_{T^ix}(T^{i+j+\ell-n_k}A_k \cap
{R}_k)=\sigma_x(T^{j+\ell-n_k}A_k\cap T^{-i}{R}_k).
\end{equation}
Also,
\begin{equation}
\label{eq:cjx}
c_j(x)=\sigma_x(T^{j+\ell}A_k \cap {R}_k).
\end{equation}
Now because $\tilde{R}_k\subset \bigcap_{i=-{n_k}}^{n_k}T^i{R}_k $, if $z \in T^{i+j+\ell-n_k} A_k \cap \tilde{R}_k$,  then,
$z \in T^{j+\ell-n_k}A_k\cap T^{-i}{R}_k$, and $z \in
T^{j+\ell}A_k \cap {R}_k$. Therefore, the symmetric difference
between $T^{j+\ell-n_k}A_k\cap T^{-i}{R}_k$ and $T^{j+\ell}A_k
\cap {R}_k$ is contained in the union of $T^{i+j+\ell-n_k} A_k \cap
\tilde{R}_k^c$ and $T^{j+\ell} A_k \cap \tilde{R}_k^c$. 
Thus, in view of (\ref{eq:cjTix}), and (\ref{eq:cjx}), 
\begin{displaymath}
|c_j(x)-c_j(T^ix)|\leq
\sigma_x(T^{j+\ell+i-n_k}A_k \cap \tilde{R}_k^c)+
\sigma_x(T^{j+\ell}A_k \cap \tilde{R}_k^c).
\end{displaymath}
 The last case, where $j+\ell \geq n_k$ and $0\leq i+j+\ell<n_k$ gives analogous bounds. So we bound $\sum_{i=0}^{n_k-1} |c_j(x)-c_j(T^ix)|$ by $2\sum_{i=0}^{n_k-1}\sigma_x(T^iA_k\cap \tilde{R}_k^c)\leq 2\sigma_x( \tilde{R}_k^c)$ and obtain the lemma.
\end{proof}

Let $d_{\mathcal{M}([0,1])}$ denote the Kantorovich-Rubinstein metric on
measures. That is
\begin{displaymath}
d_{\mathcal{M}([0,1])}(\mu,\nu)=\sup \left\{\left|\int fd\mu-\int f d\nu\right|:f \text{ is 1-Lipschitz}\right\}.
\end{displaymath}
Note, restricted to measures with total variation at most 1 it defines the same topology as the weak-* topology (on this set). 
The next lemma is an immediate consequence of the definition of $d_{\mathcal{M}([0,1])}$.
\begin{lem}\label{lemma:kr est}If $f$ is 1-Lipshitz and
$d_{\mathcal{M}([0,1])}(\sigma_x,\sigma_y)<\epsilon$ then $|A_\sigma f(x)-A_\sigma
f(y)|<\epsilon$.
\end{lem}

We say $0 \le j<n_k$ is \emph{k-good} if there exists $y_j$ in $T^jA_k$ so
that 
%\begin{itemize}
%\item $\sigma_{y_j}(T^jA_k \setminus B_j^{(k)})<\epsilon \$ 
%\item $\sum_{i:T^i(y_j)\notin B_{a_i}^{(k)}}\sigma_{y_j}(T^{a_i}A_k)<\epsilon$
at least $1-\epsilon$ proportion of the points in $T^jA_k$ have their
disintegration is $\epsilon$ close to $y_j$. 
That is
\begin{displaymath}
\lambda(\{x \in T^jA_k:
  d_{\mathcal{M}([0,1])}(\sigma_x,\sigma_{y_j})<\epsilon\}) \ge (1-\epsilon) \lambda(A_k).
\end{displaymath}
%\end{itemize} 

\begin{lem}
\label{lemma:most:good}
For all $\epsilon>0$ there exists
$k_0$ so that for all $k>k_0$ we have
$$|\{0 \le j < n_k : j \text{ is k-good }\}|>(1-\epsilon)n_k.$$
\end{lem}
\begin{proof} 
%Let  $g_k(x)=\sigma_x(\cR_k \setminus \cup_{i=0}^{n_k-1}B_i^{(k)}$. Let  $a_i(x) \in [0,n_k-1]$  
 %$j+i \equiv a_i(x) \, (mod\, n_k)$ where $x \in T^j(A_k)$ for $0\leq j<n_k$. Define 
%$$h_k(y)=\sum_{i:T^ix\notin a_i(x)} \sigma_y(T\sigma_y(T^iA_k).$$
%By \eqref{cond:bunched} we have $\int h_k(y)d\lambda, \, \int g_k d\lambda \to 0$. So by Markov's inequality (and the fact that $g_k.h_k$ are non-negative) there exists $k_{-2}$ so that $\lambda(\{y:g_k(y)<\epsilon \text{ and } h_k(y)<\epsilon\})>1-\frac \epsilon 4.$

By Lusin's Theorem there exists a compact set $K$ of measure at least
$1-\frac {\epsilon^2} 8$ so that the map $y \to \sigma_y$ is
continuous with respect to the usual metric on $[0,1]$ and the metric
$d_{\mathcal{M}([0,1])}$ on measures. Because $K$ is compact this map is uniformly
continuous and so there exists $\delta>0$ so that  $x,y \in K$ and
$|x-y|<\delta$ then $d_{\mathcal{M}([0,1])}(\sigma_x,\sigma_y)<\epsilon$. We choose
$k_0$ so that for $k>k_0$ there are $\hat{B}_i^{(k)}$ with $\diam(\hat{B}_i^{(k)})<\delta$ and $\lambda\big([0,1] \setminus
\cup_{i=0}^{n_k-1}(T^iA_k \cap \hat{B}_i^{(k)})\big)<\frac {\epsilon^2}8$ and $\lambda([0,1]\setminus \mathcal{R}_k)<\frac {\epsilon^2} 4$.  (We can do this by Condition \eqref{cond:big tower} and \eqref{cond:bunched} of rigid rank 1.) 
Let
\begin{displaymath}
\eta = \frac{1}{n_k}|\{0\leq j<n_k: \lambda\Big(T^j A_k\cap \big(K^c \cup (B_j^{(k)})^c\big)\Big)> \epsilon
\lambda(A_k)\}|. 
\end{displaymath}
Then, because the $T^j A_k$ are disjoint and of equal size and
$\bigcup_{j=0}^{n_k-1} T^j A_k = \cR_k$, it is clear that  
\begin{displaymath}
\eta \epsilon \le \frac{\lambda(K^c \cup( \cup_{i=0}^{n_k-1}T^iA_k\setminus \hat{B}_i^{(k)} )^c\cap \cR_k)}{\lambda(\cR_k)} \le \frac{\frac{\epsilon^2}4}{1-\frac {\epsilon^2}4}<\frac{\epsilon^2}2,
\end{displaymath}
and thus $\eta < \epsilon/2$. This completes the proof of the lemma.
\end{proof}

\noindent
\textbf{Notation.} 
Let 
\begin{multline*}V_j=\{x\in T^jA_k: \sigma_x\big(\cup_{i=0}^{n_k-1}(T^iA_k\setminus B_i^{(k)})\big)>(1-\epsilon)\lambda(A_k), \, \sigma_x(\tilde{R}_k^c)<\epsilon\\
 \text{ and } 
\sum_{i:T^ix\notin B_{a_i}^{(k)}} \sigma_x(T^iA_k)<\epsilon
\},
\end{multline*}
where $0\leq a_i<n_k$ is as in Lemma \ref{lemma:A close}.
If $j$ is $k$-good let 
\begin{multline*}
G_j=\{x \in T^jA_k: \lambda(\{y \in
T^jA_k:d_{\mathcal{M}(Y)}(\sigma_x,\sigma_y)<2\epsilon\})>(1-2\epsilon)\lambda(A_k)
\end{multline*}
That is, $V_j$ is a subset of $T^j(A_k)$ where  Lemma \ref{lemma:A close} gives a strong estimate and $G_j$ is the subset of $T^jA_k$ 
that are almost continuity points of the map $x \to \sigma_x$
(restricted to $T^jA_k$). We set $G_j=\emptyset$ if $j$ is not $k$-good. 

\begin{lem}
\label{lemma:our point}
For all $\epsilon>0$ there exists $k_1$ so that for all $k>k_1$ there
exists $0\leq \ell<n_k$ and $y_k\in V_\ell$ so
that 
\begin{equation}
\label{eq:good:yk}
|\{-\ell\leq j<n_k-\ell:T^j y_k \in G_{\ell+j} \cap V_{\ell+j} \text{ and }j \text{ is $k$-good}\}|>(1-13\sqrt{\epsilon})n_k. 
\end{equation}
\end{lem}
\begin{proof} If $j$ is $k$-good then 
$$\lambda(G_j)>(1-\epsilon)\lambda(A_k).$$
Let $\cR_k^* = \bigcup_{j=0}^{n_k-1} G_j$. 
Notice that  $\underset{k\to \infty}{\lim} \, \lambda(\cup_{i=0}^{n_k-1}T^iA_k)=\underset{k\to\infty}{\lim}\, \lambda(\cR_k)=1$  and so for all large enough $k$ (so that $\lambda(\cR_k)$ is close to 1 and Lemma~\ref{lemma:most:good} holds) we have 
$$\lambda(\cR_k^*)\geq  (1-\epsilon)^2\lambda(\cR_k)>1-3\epsilon.$$
By a straightforward $L^1$ estimate, we have 
\begin{multline*}
\sum_{\ell=0}^{n_k-1}\lambda(\{y\in T^\ell A_k:|\{-\ell \leq
j<n_k-\ell:G_j = \emptyset \text{ or } T^jy\not\in G_{j+\ell}\}|\geq 12\sqrt{\epsilon} n_k\} )<
%4 \frac 1 {\sqrt{\epsilon}} \lambda(G_j)
\frac{ 3 \sqrt{\epsilon}}{12}=\frac   {\sqrt{\epsilon}} 4
%< 2\lambda(\cR_k^*)<6 \epsilon.
\end{multline*}
 for all large enough $k$. 
 
 Now for the bound on $V_j$. Let $f_k(x)=\sigma_x(\tilde{\cR}_k^c)$. Let $g_k(x)=\sigma_x\big(\cup_{i=0}^{n_k-1}(T^iA_k\setminus B_i^{(k)})\big)$. Let  $a_i(x) \in [0,n_k-1]$  
 $j+i \equiv a_i(x) \, (mod\, n_k)$ where $x \in T^j(A_k)$ for $0\leq j<n_k$. Define 
$$h_k(x)=\sum_{i:T^ix\notin B^{(k)}_{a_i(x)}} \sigma_x(T^iA_k).$$
By \eqref{cond:big tower} we have $\int f_k d\lambda \to 0$ and  \eqref{cond:bunched} we have $\int h_kd\lambda, \, \int g_k d\lambda \to 0$. So by a straightforward $L^1$ estimate (and the fact that $f_k,\, g_k,\, h_k$ are non-negative)  
\begin{equation}\label{eq:V big}\lambda(\{y:f_k(y)<\epsilon, \, g_k(y)<\epsilon \text{ and } h_k(y)<\epsilon\})>1-\frac \epsilon 4
\end{equation} for all large $k$.  
%By \ref{cond:bunched} of rigid rank 1 we have that
%$$\underset{k \to \infty}{\lim}\, \lambda(\{y: \sigma_y(\cup_{i=0}^{n_k-1}T^iA_k \setminus B_i^{(k)})>\frac \epsilon 4 \}) =0$$
Therefore $ \sum_{\ell=0}^{n_k-1}\lambda(\{y\in T^\ell A_k:|\{-\ell \leq
j<n_k-\ell: T^jy\not\in V_{j+\ell}\}|\geq \sqrt{\epsilon} n_k\} )<
\frac{  \sqrt{\epsilon}}{4}.$

%Recalling that by 
%Lemma~\ref{lemma:srank est} we have $\underset{k \to \infty}{\lim}\,
%\lambda(\tilde{\cR}_k^c)= 0$ and so for $k$ large enough,
%$$\lambda(\{y: \sigma_y(\tilde{\cR}_k)>\epsilon\})<\frac 1 3.$$
%Thus, we can pick $y_k$ satisfying the conditions of the lemma. 
\end{proof}

\begin{proof}[Proof of Theorem \ref{theorem:SOT close}] 
For each $k$ large enough so that Lemmas \ref{lemma:most:good} and \ref{lemma:our point} hold  and $\lambda(\cR_k^c)<\epsilon$, let  $y_{k}$ be as in the statement of Lemma~\ref{lemma:our point} and in particular, it is in $T^\ell A_k$ for some $0\leq \ell<n_k$.
\medskip

\noindent
\textit{Step 1:} We show that for all 1-Lipschitz functions $f$ with $\|f\|_{\sup}\leq 1$ we have 
$$\underset{k \to \infty}{\lim} \, \|A_\sigma f-\sum_{i=0}^{n_k-1}c_i(y_{k})U_T^if\|_2=0.$$
First, observe that because $T^jy_k \in V_{j+\ell}$, for some $\ell$ and $j$, Lemma~\ref{lemma:A close} and the fact that $\|f\|_{\sup}\leq 1$ imply, 
\begin{equation}\label{eq:compare to coeff}
|A_\sigma f(T^jy_k)-\sum_{i=0}^{n_k-1}c_i(T^jy_k)f(T^{i+j}y_k)|<
4\epsilon.
\end{equation}
% By our assumption that $\sigma_{y_k}(\cup_{i=0}^{n_k-1}T^iA_k \setminus B_i^{(k)})$ and  $\sigma_{y_k}(\tilde{\cR}_k^c)<\epsilon$ we have
% $$ |A_\sigma f(T^jy_k)-\sum_{i=0}^{n_k-1}c_i(T^jy_k)f(T^{i+j}y_k)|<4\epsilon.$$
  From Lemma~\ref{lemma:kr est} we 
have that if $x$ satisfies 
\begin{equation}\label{eq:kr close}d_{\mathcal{M}(Y)}(\sigma_x,\sigma_{T^jy_k})<2\epsilon \end{equation}
 then 
\begin{multline*}|A_\sigma f(x)-\sum_{i=0}^{n_k-1}c_i(T^jy_k)f(T^{i+j}y_k)|\leq |A_{\sigma}f(x)-A_{\sigma}f(y_k)|+\\
 |A_\sigma f(T^jy_k)-\sum_{i=0}^{n_k-1}c_i(T^jy_k)f(T^{i+j}y_k)|<2\epsilon+4\epsilon= 6\epsilon.
\end{multline*}
%Let $V=\cup_{j=0}^{n_k-1}V_j$. % denote the set of $x$ 
%\begin{itemize}
%\item satisfying (\ref{eq:kr close}),
%\item with 
%$x \in T^{\ell+j}I_k \cap \hat{\cR}_k$ for $-\ell \leq
%j<n_k-\ell$ and 
%\item $\sigma_x(\cup_{i=0}^{n_k-1}T^iA_k\setminus B_i^{(k)})<\epsilon.$
%\end{itemize}
%Then, for $x \in V$,  
%$T^{i}x,T^{i+j}y_k \in T^{i+\ell+j\, (mod \, n_k)}I_k$ for all $0\leq i<n_k$ since $-n_k<i, \, i+j<n_k$ (by \eqref{eq:property:hatRk}). 
For any $x$  satisfying \eqref{eq:kr close},
\begin{multline*}
|A_\sigma f(x)-\sum_{i=0}^{n_k-1}c_i(T^jy_k)f(T^ix)|\leq |A_\sigma f(x)-\sum_{i=0}^{n_k-1}c_i(T^jy_k)f(T^{i+j}y_k)|+\\
 |\sum_{i=0}^{n_k-1}c_i(T^jy_k)f(T^{i+j}y_k)-\sum_{i=0}^{n_k-1}c_i(T^jy_k)f(T^ix)|<
6\epsilon+\\
\epsilon+ \sum_{i:T^ix \notin B_{a_i(x)}^{(k)}}\sigma_x(T^iA_k)+\sum_{i:T^{i+j}y_k \notin B_{a_i(y_k)}^{(k)}}\sigma_{T^{j}(y_k)}(T^iA_k).
\end{multline*}
%Note in the second inequality we use \eqref{eq:compare to coeff}. 
Now if $x,T^j y_k\in V_{\ell+j}$ we have that this is at most $9 \epsilon$. Let 
$$\hat{V}=\cup_{j\in [-\ell,n_k-\ell):T^j y_k \in V_{\ell+j}}V_{\ell+j}.$$%because $x \in V$ and there exists $\ell$ with $|\ell-j|<13\sqrt{\epsilon}$ so that $T^{\ell}y_k \in V$ we have that this is at most 
%$7\epsilon +\epsilon+13\sqrt{\epsilon}+\epsilon$. \marginpar{\jcc{More detail?}}

%By our condition that $y_k$ satisfies Lemma \ref{lemma:our point} and that $x \in V$ this is at most $8\epsilon$. 
%Recalling that by assumption $\diam(I_k)<\epsilon$ and 
By Lemma
\ref{lemma:other indices}  we have 
  $$\int_{\hat{V}}|A_\sigma f(x)-\sum_{j=0}^{n_k-1}c_j(y_k)f(T^jx)|d\lambda(x)\leq 9 \epsilon+ 2\sigma_{y_k}(\tilde{R}_k^c)\leq
   %(9\epsilon+13\sqrt{\epsilon})^2.
9 \epsilon+2\epsilon.$$
Since $|A_{\sigma}f(x)-\sum_{j=0}^{n_k-1}c_j(y_k)f(T^jx)|\leq 2$ for all $x$, by H\"older's inequality  
$$\int_{\hat{V}}|A_\sigma f(x)-\sum_{j=0}^{n_k-1}c_j(y_k)f(T^jx)|^2d\lambda(x)\leq
   %(9\epsilon+13\sqrt{\epsilon})^2.
2\cdot 11\epsilon.$$

    Since $y_k$ satisfies the assumptions of Lemma~\ref{lemma:our point}, %and \eqref{eq:V big}
     we have that 
\begin{equation}\label{eq:V big}\lambda(\{z: z\notin V_{y_k} \text{ or } z \text{ does not satisfy } \eqref{eq:kr close})<13\sqrt{\epsilon} n_k \lambda(A_k)+\epsilon+\lambda(\cR^c).
\end{equation}
  Estimating trivially on $V^c$ we have
\begin{multline*}
\|A_\sigma f-\sum_{j=0}^{n_k-1}c_j(y_k)f\circ T^j\|_2^2=\int_0^1
  |A_\sigma f(x)-\sum_{j=0}^{n_k-1}c_j(y_k)f(T^jx)|^2 \, d\lambda(x)
  \leq \\ \leq
  2 \cdot \big(13\epsilon+13\sqrt{\epsilon}\big).
\end{multline*}
Since $\|f\|_{\sup}\leq 1$ and $\epsilon$ is arbitrary this establishes Step 1.
\medskip

\noindent  
\textit{Step 2:} Completing the proof. 

Step 1 establishes pointwise convergence for a subset of $L^2$ with dense span. Because the linear operators in our sequence have uniformly bounded $L^2$ operator norm (in fact bounded by 1) this gives pointwise convergence on all of $L^2$; that is, SOT convergence.
%\marginpar{\jcc{Can we just write ``The proof follows by linearity and that Step 1 give the limit on subset of $L^2$ with dense span."}}

%The idea of the proof is that by step 1 and linearity we have the limit on a dense set in $L^2$. Since the functions on $L^2$ we consider have operator norm uniformly bounded (by 1) they are an equicontinuous family and so convergence on a dense set implies convergence. 

%To complete the formal proof of the theorem,  observe that for any $z$
%we have $\sum c_i(z)=\sum |c_i(z)|\leq \sigma_z([0,1])$ and we may
%assume that $\sigma_z([0,1])= 1$.\footnote{It is 1 for all but a
%  measure zero
%  set of $z$ and we may change the disintegration on this zero set.}
%So
%\begin{displaymath}
%\left\|\sum_{i=0}^{n_k-1}c_i(y_k)U_T^i\right\|_{op}\leq 1 \qquad\text{for all $k$}.
%\end{displaymath}
%Therefore since we have shown $\underset{k \to \infty}{\lim} \,
%\|A_\sigma f- \sum_{i=0}^{n_k-1}c_i(y_k)U_T^if\|_2 =0$ for a set of
%$f$ with dense span in $L^2$ (that is 1-Lipschitz functions with
%$\|f\|_{\sup}\leq 1$), we know that for all $f \in L^2 $  we have that
%$\underset{k \to \infty}{\lim} \, \|A_\sigma f- \sum_{i=0}^{n_k-1}c_i(y_k)U_T^if\|_2 =0.$ This is the definition of strong operator convergence. 
\end{proof}

\begin{proof}[Proof of Corollary \ref{cor:WOT close}]
Let $\hat{\delta}_p$ denote the point mass at $p$.
By the proof of
the theorem that there exists $z$ (it is $T^jy_k$ in the proof) so that
\begin{displaymath}
d_{\mathcal{M}(Y)}(\sigma_x,\sum_{j=0}^{n_k-1}c_j(z)\hat{\delta}_{(x,T^ix)})<8\epsilon
\end{displaymath}
for all $x \in V$.   By
(\ref{eq:V big}) we may assume $\lambda(V^c)$ is as small as we want. The corollary follows.
\end{proof}

\end{document}